\input miniltx \input graphicx.sty \input amssym \font \bbfive = bbm5 \font \bbeight = bbm8 \font \bbten = bbm10 \font \rs = rsfs10 \font \rssmall =
rsfs10 scaled 833 \font \eightbf = cmbx8 \font \eighti = cmmi8 \skewchar \eighti = '177 \font \fouri = cmmi5 scaled 800 \font \eightit = cmti8 \font
\eightrm = cmr8 \font \eightsl = cmsl8 \font \eightsy = cmsy8 \skewchar \eightsy = '60 \font \eighttt = cmtt8 \hyphenchar \eighttt = -1   \font \sixi = cmmi6 \skewchar \sixi = '177 \font \sixrm = cmr6 \font \sixsy = cmsy6 \skewchar \sixsy = '60 \font \tensc
= cmcsc10   \scriptfont \bffam = \bbeight \scriptscriptfont \bffam = \bbfive \textfont \bffam = \bbten
\newskip \ttglue \def \eightpoint {\def \rm {\fam 0 \eightrm }\relax \textfont 0 = \eightrm \scriptfont 0 = \sixrm \scriptscriptfont 0 = \fiverm
\textfont 1 = \eighti \scriptfont 1 = \sixi \scriptscriptfont 1 = \fouri \textfont 2 = \eightsy \scriptfont 2 = \sixsy \scriptscriptfont 2 = \fivesy
\textfont 3 = \tenex \scriptfont 3 = \tenex \scriptscriptfont 3 = \tenex \def \it {\fam \itfam \eightit }\relax \textfont \itfam = \eightit \def \sl
{\fam \slfam \eightsl }\relax \textfont \slfam = \eightsl \def \bf {\fam \bffam \eightbf }\relax \textfont \bffam = \bbeight \scriptfont \bffam =
\bbfive \scriptscriptfont \bffam = \bbfive \def \tt {\fam \ttfam \eighttt }\relax \textfont \ttfam = \eighttt \tt \ttglue = .5em plus.25em minus.15em
\normalbaselineskip = 9pt \def \MF {{\manual opqr}\-{\manual stuq}}\relax \let \sc = \sixrm \let \big = \eightbig \let \rs = \rssmall \setbox \strutbox
= \hbox {\vrule height7pt depth2pt width0pt}\relax \normalbaselines \rm } \def \setfont #1{\font \auxfont =#1 \auxfont } \def \withfont #1#2{{\setfont
{#1}#2}}    \def \TRUE {Y} \def \FALSE {N} \def \ifundef
#1{\expandafter \ifx \csname #1\endcsname \relax } \def \undefrule {\kern 2pt \vrule width 2pt height 5pt depth 0pt \kern 2pt} \def \UndefLabels {}
\def \possundef #1{\ifundef {#1}\undefrule {\eighttt #1}\undefrule \global \edef \UndefLabels {\UndefLabels #1\par } \else \csname #1\endcsname \fi }
\newcount \secno \secno = 0 \newcount \stno \stno = 0 \newcount \eqcntr \eqcntr = 0 \ifundef {showlabel} \global \def \showlabel {\FALSE } \fi \ifundef
{auxwrite} \global \def \auxwrite {\TRUE } \fi \ifundef {auxread} \global \def \auxread {\TRUE } \fi \def \define #1#2{\global \expandafter \edef \csname
#1\endcsname {#2}} \long \def \error #1{\medskip \noindent {\bf ******* #1}} \def \fatal #1{\error {#1\par Exiting...}\end } \def \advseqnumbering
{\global \advance \stno by 1 \global \eqcntr =0} \def \current {\ifnum \secno = 0 \number \stno \else \number \secno \ifnum \stno = 0 \else .\number
\stno \fi \fi } \def \rem #1{\vadjust {\vbox to 0pt{\vss \hfill \raise 3.5pt \hbox to 0pt{ #1\hss }}}} \font \tiny = cmr6 scaled 800 \def \deflabel
#1#2{\relax \if \TRUE \showlabel \rem {\tiny #1}\fi \ifundef {#1PrimarilyDefined}\relax \define {#1}{#2}\relax \define {#1PrimarilyDefined}{#2}\relax
\if \TRUE \auxwrite \immediate \write 1 {\string \newlabe l {#1}{#2}}\fi \else \edef \old {\csname #1\endcsname }\relax \edef \new {#2}\relax \if \old
\new \else \fatal {Duplicate definition for label ``{\tt #1}'', already defined as ``{\tt \old }''.}\fi \fi } 
\def \label #1 {\deflabel {#1}{\current }} \def \equationmark #1 {\ifundef {InsideBlock} \advseqnumbering \eqno {(\current )} \deflabel {#1}{\current
} \else \global \advance \eqcntr by 1 \edef \subeqmarkaux {\current .\number \eqcntr } \eqno {(\subeqmarkaux )} \deflabel {#1}{\subeqmarkaux } \fi \if
\TRUE \showlabel \hbox {\tiny #1}\fi } \def \lbldeq #1 $$#2$${\ifundef {InsideBlock}\advseqnumbering \edef \lbl {\current }\else \global \advance \eqcntr
by 1 \edef \lbl {\current .\number \eqcntr }\fi $$ #2 \deflabel {#1}{\lbl }\eqno {(\lbl )} $$} \def \split #1.#2.#3.#4;{\global \def \parone {#1}\global
\def \partwo {#2}\global \def \parthree {#3}\global \def \parfour {#4}} \def \NA {NA} \def \ref #1{\split #1.NA.NA.NA;(\possundef {\parone }\ifx \partwo
\NA \else .\partwo \fi )}  \newcount \bibno \bibno = 0  \def \Bibitem #1 #2; #3;
#4 \par {\smallbreak \global \advance \bibno by 1 \item {[\possundef {#1}]} #2, {``#3''}, #4.\par \ifundef {#1PrimarilyDefined}\else \fatal {Duplicate
definition for bibliography item ``{\tt #1}'', already defined in ``{\tt [\csname #1\endcsname ]}''.} \fi \ifundef {#1}\else \edef \prevNum {\csname
#1\endcsname } \ifnum \bibno =\prevNum \else \error {Mismatch bibliography item ``{\tt #1}'', defined earlier (in aux file ?) as ``{\tt \prevNum }''
but should be ``{\tt \number \bibno }''.  Running again should fix this.}  \fi \fi \define {#1PrimarilyDefined}{#2}\relax \if \TRUE \auxwrite \immediate
\write 1 {\string \newbi b {#1}{\number \bibno }}\fi } \def \jrn #1, #2 (#3), #4-#5;{{\sl #1}, {\bf #2} (#3), #4--#5} \def \Article #1 #2; #3; #4 \par
{\Bibitem #1 #2; #3; \jrn #4; \par } \def \references {\begingroup \bigbreak \eightpoint \centerline {\tensc References} \nobreak \medskip \frenchspacing }
\catcode `\@=11 \def \citetrk #1{{\bf \possundef {#1}}} \def \c@ite #1{{\rm [\citetrk {#1}]}} \def \sc@ite [#1]#2{{\rm [\citetrk {#2}\hskip 0.7pt:\hskip
2pt #1]}} \def \du@lcite {\if \pe@k [\expandafter \sc@ite \else \expandafter \c@ite \fi } \def \cite {\futurelet \pe@k \du@lcite } \catcode `\@=12 \def
\Headlines #1#2{\nopagenumbers \headline {\ifnum \pageno = 1 \hfil \else \ifodd \pageno \tensc \hfil \lcase {#1} \hfil \folio \else \tensc \folio \hfil
\lcase {#2} \hfil \fi \fi }} \def \title #1{\medskip \centerline {\withfont {cmbx12}{\ucase {#1}}}}  \long \def \Quote #1\endQuote {\begingroup \leftskip 35pt \rightskip 35pt \parindent 17pt
\eightpoint #1\par \endgroup } \long \def \Abstract #1\endAbstract {\vskip 1cm \Quote \noindent #1\endQuote }   \def \Note #1{\footnote {}{\eightpoint #1}} \def \Date #1 {\Note
{\it Date: #1.}} \newcount \auxone \newcount \auxtwo \newcount \auxthree \def \currenttime {\auxone =\time \auxtwo =\time \divide \auxone by 60 \auxthree
=\auxone \multiply \auxthree by 60 \advance \auxtwo by -\auxthree \ifnum \auxone <10 0\fi \number \auxone :\ifnum \auxtwo <10 0\fi \number \auxtwo } \def
\today {\ifcase \month \or January\or February\or March\or April\or May\or June\or July\or August\or September\or October\or November\or December\fi
{ }\number \day , \number \year }  \def \hojeExtenso {\number \day \ de \ifcase \month \or
janeiro\or fevereiro\or mar\c co\or abril\or maio\or junho\or julho\or agosto\or setembro\or outubro\or novembro\or decembro\fi \ de \number \year }  \def \part #1#2{\vfill \eject \null \vskip 0.3\vsize \withfont {cmbx10 scaled 1440}{\centerline {PART #1}
\vskip 1.5cm \centerline {#2}} \vfill \eject }   \def \fix {\smallskip \noindent $\blacktriangleright $\kern 12pt} 
 \def \ucase #1{\edef \auxvar {\uppercase {#1}}\auxvar } \def \lcase #1{\edef \auxvar {\lowercase {#1}}\auxvar } \def
\emph #1{{\it #1}\/} \def \section #1 \par {\global \advance \secno by 1 \stno = 0 \goodbreak \bigbreak \noindent {\bf \number \secno .\enspace #1.}
\nobreak \medskip \noindent } \def \state #1 #2\par {\begingroup \def \InsideBlock {} \medbreak \noindent \advseqnumbering {\bf \current .\enspace
#1.\enspace \sl #2\par }\medbreak \endgroup } \def \definition #1\par {\state Definition \rm #1\par } \newcount \CloseProofFlag   \long \def \Proof #1\endProof {\begingroup \def
\InsideBlock {} \global \CloseProofFlag =0 \medbreak \noindent {\it Proof.\enspace }#1 \ifnum \CloseProofFlag =0 \hfill $\endproofmarker $ \looseness =
-1 \fi \medbreak \endgroup }  \def \explica #1#2{\mathrel {\buildrel \hbox {\sixrm #1} \over #2}} \def \explain #1#2{\explica
{\ref {#1}}{#2}}  \def \=#1{\explain {#1}{=}} \def \pilar #1{\vrule height #1 width 0pt}
 \newcount \fnctr \fnctr = 0 \def \fn #1{\global \advance \fnctr by 1 \edef \footnumb {$^{\number \fnctr
}$}\relax \footnote {\footnumb }{\eightpoint #1\par \vskip -10pt}} \def \text #1{\hbox {#1}}   \def
\item #1{\par \noindent \kern 1.1truecm\hangindent 1.1truecm \llap {#1\enspace }\ignorespaces }  \def \Item #1{\smallskip \item
{{\rm #1}}} \newcount \zitemno \zitemno = 0 \def \izitem {\global \zitemno = 0} 
 \def \zitemplus {\global \advance \zitemno by 1 \relax } \def \rzitem {\romannumeral \zitemno } \def \rzitemplus {\zitemplus
\rzitem } \def \zitem {\Item {{\rm (\rzitemplus )}}}   \newcount \nitemno \nitemno = 0  \def \nitem {\global \advance \nitemno by 1 \Item {{\rm (\number \nitemno )}}} \newcount \aitemno \aitemno = -1 \def \boxlet #1{\hbox to
6.5pt{\hfill #1\hfill }}  \def \aitemconv {\ifcase \aitemno a\or b\or c\or d\or e\or f\or g\or h\or i\or j\or k\or l\or m\or
n\or o\or p\or q\or r\or s\or t\or u\or v\or w\or x\or y\or z\else zzz\fi } \def \aitem {\global \advance \aitemno by 1\Item {(\boxlet \aitemconv )}}
 \def \deflabeloc #1#2{\deflabel {#1}{\current .#2}{\def \showlabel {\FALSE }\deflabel {Local#1}{#2}}} \def \lbldzitem
#1 {\zitem \deflabeloc {#1}{\rzitem }} \def \lbldaitem #1 {\aitem \deflabeloc {#1}{\aitemconv }} \def \aitemmark #1 {\deflabel {#1}{\aitemconv }}
\def \iItemmark #1 {\zitemmark {#1} } \def \zitemmark #1 {\deflabel {#1}{\current .\rzitem }{\def \showlabel {\FALSE }\deflabel {Local#1}{\rzitem
}}} \def \Case #1:{\medskip \noindent {\tensc Case #1:}} \def \<{\left \langle \vrule width 0pt depth 0pt height 8pt } \def \>{\right \rangle }
\def \({\big (} \def \){\big )}  \def \and {\hbox {\quad and \quad }} 
  \def \IFF {\kern 7pt\Leftrightarrow \kern 7pt} \def \IMPLY {\kern
7pt \Rightarrow \kern 7pt} \def \for #1{\quad \forall \,#1} \def \endproofmarker {\square } \def \"#1{{\it #1}\/} 
 \def \*{\otimes } \def \caldef #1{\global \expandafter \edef \csname #1\endcsname {{\cal #1}}} \def \mathcal #1{{\cal #1}} \def
\bfdef #1{\global \expandafter \edef \csname #1\endcsname {{\bf #1}}} \bfdef N \bfdef Z \bfdef C \bfdef R  \def \exists
{\mathchar "0239\kern 1pt } \if \TRUE \auxread \IfFileExists {\jobname .aux}{\input \jobname .aux}{\null } \fi \if \TRUE \auxwrite \immediate \openout
1 \jobname .aux \fi \def \close {\if \empty \UndefLabels \else \message {*** There were undefined labels ***} \medskip \noindent ****************** \
Undefined Labels: \tt \par \UndefLabels \fi \if \TRUE \auxwrite \closeout 1 \fi \par \vfill \supereject \end } \def \Caixa #1{\setbox 1=\hbox {$#1$\kern
1pt}\global \edef \tamcaixa {\the \wd 1}\box 1} \def \caixa #1{\hbox to \tamcaixa {$#1$\hfil }}   \def \medcup {\mathop {\mathchoice {\raise 1pt \hbox {$\mathchar "1353$}}{\mathchar "1353}{\mathchar "1353}{\mathchar "1353}}}
\def \medcap {\mathop {\mathchoice {\raise 1pt \hbox {$\mathchar "1354$}}{\mathchar "1354}{\mathchar "1354}{\mathchar "1354}}}         \def \cl #1 #2 #3 {#1, & \hbox {#2 } #3\hfill \crr } \def \paper #1#2#3{ \hsize #1truemm
\advance \hsize by -#3truemm  \advance \hsize by -#3truemm \vsize #2truemm   \advance \vsize by -#3truemm  \advance \vsize by -#3truemm \hoffset =-1truein
\advance \hoffset by #3truemm \voffset =-1truein \advance \voffset by #3truemm } \def \qt $#1${``$#1$''} \def \qt $#1${$#1$} \def \nd /{non-degenerate}
\def \Nd /{Non-degenerate} \def \nega {\neg \,} \def \emptyset {\varnothing } \def \cf [#1]#2{({\it cf}.~\cite [#1]{#2})} \def \ideal #1{\big \langle
#1\big \rangle } \newcount \currentgrouplevel

\centerline {\bf TIGHT AND COVER-TO-JOIN REPRESENTATIONS OF} \smallskip \centerline {\bf SEMILATTICES AND INVERSE SEMIGROUPS}

\bigskip \centerline {\tensc R. Exel\footnote {$^{\ast }$}{\eightrm Universidade Federal de Santa Catarina and University of Nebraska at Lincoln.} }

\Abstract We discuss the relationship between tight and cover-to-join representations of semilattices and inverse semigroups, showing that a slight
extension of the former, together with an appropriate selection of co-domains, makes the two notions equivalent.  As a consequence, when constructing
universal objects based on them, one is allowed to substitute cover-to-join for tight and vice-versa.  \endAbstract

\footnote {}{\eightrm Date: March 7, 2019.}

\section Introduction

Exactly twelve years ago, to be precise on March 7,  2007,  I posted a paper on the arXiv \cite {actionspre} describing the notion of \emph {tight
representations} of semilattices and inverse semigroups, which turned out to have many applications and in particular proved to be useful to give a
unified perspective to a significant number of  C*-algebras containing a preferred generating set of partial isometries (\cite {Boava}, \cite {DonMil},
\cite {actions}, \cite {EGS}, \cite {EP}, \cite {ExelStar}, \cite {StarlingTwo}, \cite {StarlingThree}).

The notion of tight representations (described below for the convenience of the reader)  is slightly involving as it depends on the analysis of certain
pairs of finite sets $X$ and $Y$, but it becomes much simplified  when $X$ is a singleton and $Y$ is empty (see \cite [Proposition 11.8]{actions}).
In this simplified form it has been rediscovered and used in many subsequent works (e.g. \cite {DonMil},  \cite {Lawson},  \cite {LJ},  \cite {LV})
under the name of \emph {cover-to-join} representations.

The notion of cover-to-join representations, requiring a smaller  set of conditions, is consequently weaker and, as it turns out, strictly weaker,
than the original notion of tightness.  Nevertheless, besides being easier to  formulate,  the notion of cover-to-join representations has the
advantage of being applicable to representations taking values  in \emph {generalized} Boolean algebras, that is, Boolean algebras without a unit.
Explicitly mentioning the operation of complementation, tight representations only make sense for unital Boolean algebras.

The goal of this note is to describe an  attempt to reconcile the notions of tight and cover-to-join representations:  slightly extending the former,
and adjusting for the appropriate co-domains, we show that, after all, the two notions coincide.

One of the  main practical  consequences of this fact is that the difference between the two notions becomes irrelevant  for the purpose of constructing
universal objects based on them,  such as the completion of an inverse semigroup recently introduced in \cite {LV}.  We are moreover able to fix a slight
imprecision in the proof of \cite [Theorem 2.2]{DonMil}, at least as far as its consequence that the universal C*-algebras for tight vs.~cover-to-join
representations are isomorphic.

\section Generalized Boolean algebras

We begin by recalling the well known notion of generalized Boolean algebras.

\definition A \emph {generalized Boolean algebra} \cite [Definition 5]{Stone} is a set $B$ equipped with binary operations \qt $\wedge $ and \qt $\vee $,
and containing an element $0$,  such that for every $a$, $b$ and $c$ in $B$, one has that \def \tbox #1 / #2 / #3 /{\hbox to 2.6cm{(#1)\hfill } \hbox to
4.5cm{\hfill #2,\hfill }  \ and \  \hbox to 4.5cm{\hfill #3, \hfill }} \izitem \zitem (commutativity) \  $a\vee b = b\vee a$, and $a\wedge b = b\wedge
a$, \zitem (associativity) \ $(a\wedge b)\wedge c = a\wedge (b\wedge c)$, \zitem (distributivity) \ $a\wedge (b\vee c) = (a\wedge b)\vee (a\wedge c)$,
\zitem $a\vee 0 =a$, \zitem (relative complement) \ if $a = a\wedge b$, there is an element $x$ in $B$, such that $x \vee a = b$,  and $x \wedge a =
0$, \zitem $a\vee a =a = a\wedge a$.

\medskip It follows that (ii) and (iii) also hold with \qt $\vee $ and \qt $\wedge $ interchanged, meaning that \qt $\vee $ is associative \cite
[Theorems 55 \& 14]{Stone}, and that \qt $\vee $ distributes over \qt $\wedge $ \cite [Theorems 55 \& 11]{Stone}.

\medskip When $a=a\wedge b $, as in (v), one writes $a\leq b$.  It is then easy to see that $\leq $ is a partial order on $B$.

\medskip The element $x$ referred to in (v) is called the \emph {relative complement} of $a$ in $b$, and it is usually denoted  $b\setminus a$.

\definition \cf [Theorem 56]{Stone}  A generalized Boolean algebra $B$ is called a \emph {Boolean algebra} if there exists an element $1$ in $B$,
such that   $a\wedge 1 = a$, for every $a$ in $B$.

For Boolean algebras,  the complement of an element $a$ relative to $1$ is often denoted $\nega a$.

\medskip Recall that an \emph {ideal} of a generalized Boolean algebra $B$ is any nonempty subset  $C$ of $B$ which is  closed under \qt $\vee $, and
such that $$ a\leq b \in C \IMPLY a\in C.  $$ Such an ideal is evidently also closed under \qt $\wedge $ and under relative complements, so it is a
generalized Boolean algebra in itself.

Given any nonempty subset $S$ of $B$, notice that the subset $C$ defined by $$ C=\big \{a\in B:  \textstyle a\leq \bigvee _{z\in Z} z,  \hbox { for
some finite subset } Z\subseteq S\big \}, $$ is an ideal of $B$ and it is clearly the smallest ideal containing $S$,  so we shall call it the \emph
{ideal generated by} $S$,  and we shall denote it by  $\ideal S$.

\section Tight and cover-to-join representations of semilattices

\fix From now on let us  fix a semilattice $E$ (always assumed to have a zero element).

\definition A representation of $E$ in a generalized Boolean algebra $B$ is any map $\pi :E\to B$, such that \izitem \zitem $\pi (0)=0$, and \zitem $\pi
(x\wedge y)=\pi (x)\wedge \pi (y)$, for every $x$ and $y$ in $E$.

\medskip In order to spell out the definition of the notion of \emph {tight representations}, introduced in \cite {actions}, let $F$ be any subset
of $E$.  We then say that a given subset $Z\subseteq F$ is a \emph {cover} for $F$, if for every nonzero $x$ in $F$, there exists some $z$ in $Z$,
such that $z\wedge x \neq 0$.

Furthermore, if $X$ and $Y$ are finite subsets of $E$,  we let $$ E^{X,Y} = \{z\in E: z\leq x,\ \forall x\in X,\hbox { and } z\perp y,\ \forall y\in Y\}.  $$

\definition \label DefineTight \cf [Definition   11.6]{actions} A representation $\pi $ of $E$ in a Boolean algebra $B$ is said to be \emph {tight}
if, for any finite subsets $X$ and $Y$ of $E$, and for any finite cover $Z$ for $E^{X,Y}$, one has that \lbldeq CondForTight $$ \bigvee _{z\in Z} \pi
(z) = \bigwedge _{x\in X} \pi (x) \wedge \bigwedge _{y\in Y} \nega \pi (y).  $$

Observe that if $Y$ is empty and $X$ is a singleton, say $X=\{x\}$, then $$ E^{X,Y} =    E^{\{x\},\emptyset } = \{z\in E: z\leq x\}, $$ and if $Z$
is a cover for this set, then \ref {CondForTight} reads \lbldeq coverJoinn $$ \bigvee _{z\in Z} \pi (z) = \pi (x).  $$

To check that a given representation is tight, it is not enough to verify \ref {coverJoinn}, as it is readily  seen by considering the example in which
$E=\{0, 1\}$ and $B$ is any Boolean algebra containing an element $x\neq 1$.   Indeed, the map $\pi :E\to B$ given by $\pi (0)=0$, and $\pi (1)=x$,
satisfies all instances of \ref {coverJoinn}  even thought it is not tight.  The reader might wonder if the fact that $\pi $ fails to preserve the unit is
playing a part in this counter-example, but it is also easy to find examples of cover-to-join representations of non-unital semilattices which are not tight.

Representations $\pi $ satisfying \ref {coverJoinn}  whenever $Z$ is a cover for $E^{\{x\}, \emptyset }$ have been considered in \cite [Proposition
11.8]{actions}, and they have been called \emph {cover-to-join} representations in \cite {DonMil}.

It is a trivial matter to prove that a cover-to-join representation satisfies \ref {CondForTight} whenever $X$ is nonempty (see the proof of \cite
[Lemma 11.7]{actions}), so the question of whether a cover-to-join representation is indeed tight rests on verifying \ref {CondForTight} when $X$
is empty.  In this case, and assuming that $Z$ is a cover for $E^{\emptyset , Y}$, it is easy to see that $Z\cup Y$ is a cover for the whole of $E$.
Should we be dealing with a semilattice not admitting any finite cover, this situation will therefore never occur, that is, one will never be required
to check \ref {CondForTight} for an empty set $X$, hence every cover-to-join representation is automatically tight.

This has in fact  already been observed in \cite [Proposition 11.8]{actions}, which says that every cover-to-join representation is tight in case $E$
does not admit any finite cover, as we have just discussed, but also if $E$ contains a finite set $X$ such that \lbldeq SupCond $$ \bigvee _{x\in X}
\pi (x) = 1.  $$

The latter condition is useful for dealing with characters,  i.e.~with representations of $E$ in the Boolean algebra $\{0, 1\}$, because the requirement
that a character be nonzero immediately implies \ref {SupCond}, so again cover-to-join suffices to prove tightness.

\medskip On the other hand, an advantage of the notion of cover-to-join representations is that it makes sense for representations in generalized Boolean
algebras, while the reference to the unary operation \qt $\nega $ in \ref {CondForTight} precludes it from being applied when the target algebra lacks
a unit, that is, for a representation into a generalized Boolean algebra.

Again referring to the occurrence of \qt $\nega $ in \ref {CondForTight}, observe that if $X$ is nonempty, then the right hand side of \ref {CondForTight}
lies in the ideal of $B$ generated by the range of $\pi $.  This is because, even though $\nega \pi (y)$ is not necessarily in $\ideal {\pi (E)}$,
this term will appear besides  $\pi (x)$, for some $x$ in $X$, and hence $$ \pi (x)\wedge \nega \pi (y) =   \pi (x)\setminus \big (\pi (x)\wedge \pi
(y)\big ) \in \ideal {\pi (E)}.  $$

This means that:

\state Proposition If $E$ is a semilattice not admitting any finite cover then, whenever $X$ and $Y$ are finite subsets of $E$, and $Z$ is a finite
cover of $E^{X,Y}$, the right hand side of \ref {CondForTight} lies in $\ideal {\pi (E)}$.

As a consequence we see that definition  \ref {DefineTight} may be safely applied to a representation of $E$ in a generalized Boolean algebra, as long
as $E$ does not admit a finite cover:  despite the occurrence of \qt $\nega $ in \ref {CondForTight}, once its right hand side is expanded, it may
always be expressed  in terms of relative complements, hence avoiding the use of the missing unary operation $\nega $.

We may therefore consider the following slight generalization of the notion of tight representations:

\definition \label DefineTightTwo A representation $\pi $ of $E$ in a generalized Boolean algebra $B$ is said to be \emph {tight} if, either $B$ is
a Boolean algebra and $\pi $ is tight in the sense of \ref {DefineTight}, or the following two conditions are verified: \izitem \zitem $E$ admits no
finite cover, and \zitem \ref {CondForTight} holds for any finite subsets $X$ and $Y$ of $E$, and for any finite cover $Z$ for $E^{X,Y}$.

As already stressed, despite the occurrence of \qt $\nega $ in \ref {CondForTight}, condition \ref {DefineTightTwo.ii} will always make sense in a
generalized Boolean algebra.

So here is a result that perhaps may be used to reconcile the notions of tightness and cover-to-join representations:

\state Theorem \label MainSlat Let $\pi $ be a representation of the semilattice $E$ in the generalized Boolean algebra $B$.  Then \izitem \zitem if
$\pi $ is tight then it is also cover-to-join, \zitem if $\pi $ is cover-to-join then there exists an ideal $B'$ of $B$, containing the range of $\pi $,
such that, once $\pi $ is seen as a representation of $E$ in $B'$, one has that  $\pi $ is tight.

\Proof Point (i) being immediate, let us prove (ii).  Under the assumption that $E$ does not admit any finite cover, we have that $\pi $ is tight as
a representation into $B'=B$, by \cite [Proposition 11.8]{actions}, or rather by its obvious adaptation to generalized Boolean algebras.

It therefore remains to prove (ii) in case $E$ does admit a finite cover, say $Z$.  Setting \lbldeq DefE $$ e = \bigvee _{z\in Z} \pi (z), $$ we claim
that \lbldeq Boundd $$ \pi (x)\leq e, \for x \in E.  $$ To see this, pick $x$ in $E$ and notice that, since $Z$ is a cover for $E$, we have in particular
that the set $$ \{z\wedge x: z\in Z\} $$ is a cover for $x$, so the cover-to-join property of $\pi $ implies that $$ \pi (x) = \bigvee _{z\in Z} \pi
(z\wedge x) \leq \bigvee _{z\in Z} \pi (z) = e, $$ proving \ref {Boundd}.  We therefore let $$ B'=\{a\in B: a\leq e\}, $$ which is evidently an ideal
of $B$ containing the range of $\pi $ by \ref {Boundd}.

By \ref {DefE} we then have that $\pi $ satisfies \cite [Lemma 11.7.(i)]{actions}, as long as we see $\pi $ as a representation of $E$ in $B'$, whose
unit is clearly $e$.  The result then follows from \cite [Proposition 11.8]{actions}.  \endProof

\section \Nd / representations of semilattices

The following is perhaps the most obvious adaptation of the notion of \nd / representations extensively used in the theory of operator algebras \cite
[Definition 9.3]{Takesaki}.

\definition We shall say that a representation $\pi $ of a semilattice $E$ in a generalized Boolean algebra $B$ is \emph {\nd /} if, for every $a$ in
$B$, there is a finite subset $Z$ of $E$ such that $ a\leq \bigvee _{z\in Z} \pi (z).  $ In other words, $\pi $ is \nd / if and only if $B$ coincides
with the ideal generated by the range of $\pi $.

Observe that,  if both $E$ and $B$ have a unit, and if $\pi $ is a unital map, then $\pi $ is evidently \nd /.  More generally, if $\pi $ satisfies
\ref {SupCond}, then the same is also clearly true.

The following result says that, by adjusting the co-domain of a representation, we can always make it \nd /.

\state Proposition \label NdegNice Let $\pi $ be a representation of $E$ in the generalized Boolean algebra $B$.  Letting $C$ be the ideal of $B$
generated by the range of $\pi $, one has that $\pi $ is a \nd / representation of $E$ in $C$.

\Proof Obvious.  \endProof

For \nd / representations we have the following streamlined version of \ref {MainSlat}:

\state Corollary Let $\pi $ be a \nd / representation of the semilattice $E$ in the generalized Boolean algebra $B$.  Then $\pi $ is tight if and only
if it is cover-to-join.

\Proof The ``only if'' direction being trivial, we concentrate on the ``if'' part, so let us  assume that $\pi $ is cover-to-join.  By \ref {MainSlat}
there exists an ideal $B'$ of $B$, containing the range of $\pi $, and such that $\pi $ is tight as a representation in $B'$.  Such an ideal will therefore
contain the ideal generated by $\pi (E)$, which coincides with $B$ by hypothesis.  Therefore $B'=B$, and hence $\pi $ is tight as a representation into
its default co-domain $B$.  \endProof

\section Representations of inverse semigroups

By its very nature, the concept of a tight representation pertains to the realm of semilattices and Boolean algebras.  However, given the relevance of
the study of semilattices in the theory of inverse semigroups, tight representations have had a strong impact on the latter.

Recall that a \emph {Boolean inverse semigroup} (see \cite {TightISG} but please observe that this notion is not equivalent to the homonym studied in
\cite {Lawson} and \cite {Wehrung}) is an inverse semigroup whose idempotent semilattice $E(S)$  is a Boolean algebra.   In accordance with  what we
have been discussing up to now, it is sensible to give the following:

\definition \izitem \zitem A \emph {generalized Boolean inverse semigroup} is an inverse  semigroup   whose idempotent semilattice is a generalized
Boolean algebra.  \zitem ({\it cf}.~\cite [Definition 13.1]{actions} and \cite [Proposition 6.2]{TightISG}) If $S$ is any inverse semigroup\fn {All
inverse semigroups in this note are required to have a zero.} and $T$ is a generalized Boolean inverse semigroup, we say that a homomorphism $\pi :S\to
T$ (always assumed to preserve zero) is \emph {tight} if the restriction of $\pi $ to $E(S)$ is a tight representation into $E(T)$, in the sense of
\ref {DefineTightTwo}.  \zitem If $\pi $ is as above, we say that $\pi $ is \emph {cover-to-join} if the restriction of $\pi $ to $E(S)$ is cover-to-join.

We then have the following version of \ref {MainSlat} and  \ref {NdegNice}:

\state Corollary \label CoverTightISG Let $\pi $ be a representation of the inverse semigroup $S$ in the generalized Boolean inverse semigroup $T$.
Then \izitem \zitem if $\pi $ is tight then it is also cover-to-join, \zitem if $\pi $ is cover-to-join then there exists a generalized Boolean inverse
sub-semigroup $T'$ of $T$, containing the range of $\pi $, such that, once $\pi $ is seen as a representation of $S$ in $T'$, one has that  $\pi $
is tight.  \zitem if $\pi $ is cover-to-join,  and if the restriction of $\pi $ to $E(S)$ is \nd /,  then $\pi $ is tight.

\Proof The proof is essentially contained in the proofs of \ref {MainSlat} and \ref {NdegNice}, except maybe for the proof of (ii) under the assumption that
$E(S)$ admits  a finite cover, say $Z$.  In this case, let $e$ be as in \ref {DefE} and put $$ T'=\{t\in T: t^*t\leq e, \ tt^*\leq e\}, $$ observing that
$T'$ is clearly an inverse sub-semigroup of $T$, and that its idempotent semilattice is a Boolean algebra.  Given any $s$ in $S$, observe that  $s^*s$
lies in $E(S)$ and $$ \pi (s)^*\pi (s) =   \pi (s^*s) \leq e, $$ where the last inequality above follows as in \ref {Boundd}.  By a similar reasoning
one shows that also $\pi (s)\pi (s)^* \leq e$, so we see that $\pi (s)$ lies in $T'$, and we may then think of $\pi $ as a representation of $S$ in $T'$.
As in \ref {MainSlat}, one may now easily prove that $\pi $ becomes a tight representation into $T'$.  \endProof

\section Conclusion

As a  consequence of the above results, when defining universal objects (such as semigroups,  algebras or C*-algebras) for a class of representations
of inverse semigroups, one may safely substitute cover-to-join for tight and vice-versa.  Given the widespread use of tight representations, there are
many instances where the above principle applies.  Below we spell out one such result to concretely illustrate our point, but similar results may be
obtained as trivial reformulations of the following:

\def \ctj {C^*_{\hbox {\sixrm cover-to-join}}(S)} \def \ctj {C^*_{\buildrel {\hbox {\sixrm cover-}} \over {\hbox {\sixrm to-join}} }(S)} \def \ctight
{C^*_{\hbox {\sixrm tight}}(S)}

\state Theorem Let $S$ be an inverse semigroup and let $\ctight $ be the universal C*-algebra \cite [Theorem 13.3]{actions} for tight Hilbert space
representations of $S$ \cite [Definition 13.1]{actions}.  Also let $\ctj $ be the universal C*-algebra for cover-to-join Hilbert space representations
of $S$.  Then $$ \ctight \simeq \ctj .  $$

\Proof It suffices to prove that $\ctight $ also has the universal property for cover-to-join representations.  So let $$ \pi : S\to B(H) $$ be a
cover-to-join representation of $S$ on some Hilbert space $H$.  Should the idempotent semilattice of $S$ admit no finite covers, one has that $\pi $
is tight so there is nothing to do.  On the other hand, assuming that $Z$ is a finite cover for $E(S)$, let $e$ be as in \ref {DefE}.

Writing $H_e$ for the range of $e$ and letting $K=H_e^\perp $, we then obviously have that $H=H_e\oplus K$.    It then follows from \ref {Boundd}
that each $\pi (s)$ decomposes as a direct sum of operators $$ \pi (s) = \pi '(s)\oplus 0, $$ thus defining a representation $\pi '$ of $S$ on $H_e$
which is clearly also cover-to-join.  It is also clear that $\pi '$ is  \nd / on $E(S)$, so we have by \ref {CoverTightISG.iii} that $\pi '$ is tight.
Therefore the universal property provides a *-representation $\varphi '$ of $\ctight $ on $B(H_e)$ coinciding with $\pi '$ on the canonical image of $S$
within $\ctight $.  It then follows that $\varphi :=\varphi '\pilar {10pt}\oplus 0$  coincides with $\pi $ on $S$, concluding the proof.  \endProof

\references

\Article Boava G. Boava, G. G. de Castro, F. de L. Mortari; Inverse semigroups associated with labelled spaces and their tight spectra; Semigroup Forum,
94 (2017), no. 3, 582-609

\Article DonMil A. P. Donsig and D. Milan; Joins and covers in inverse semigroups and tight C*-algebras; Bull. Aust. Math. Soc., 90 (2014), no. 1, 121-133

\Bibitem actionspre R. Exel; Inverse semigroups and combinatorial C*-algebras (preprint version); arXiv:math/0703182v1 [math.OA], March 7, 2007

\Article actions R. Exel; Inverse semigroups and combinatorial C*-algebras; Bull. Braz. Math. Soc. (N.S.), 39 (2008), 191-313

\Article TightISG R. Exel; Tight representations of semilattices and inverse semigroups; Semigroup Forum, 79 (2009), 159-182

\Article EGS R. Exel, D. Goncalves, C. Starling; The tiling C*-algebra viewed as a tight inverse semigroup algebra; Semigroup Forum, 84 (2012), no. 2, 229-240

\Article EPdois R. Exel, E. Pardo; The tight groupoid of an inverse semigroup; Semigroup Forum, 92 (2016), no. 1, 274-303

\Article EP R. Exel, E. Pardo; Self-similar graphs, a unified treatment of Katsura and Nekrashevych $\rm C^*$-algebras; Adv. Math., 306 (2017), 1046-1129

\Article ExelStar R. Exel, C. Starling; Self-similar graph $C^*$-algebras and partial crossed products; J. Operator Theory, 75 (2016), no. 2, 299-317

\Bibitem Lawson M. V.  Lawson; Non-commutative Stone duality: inverse semigroups, topological groupoids and C ∗-algebras; Internat. J. Algebra Comput. 22
(2012), no. 6, 1250058, 47 pp

\Article LJ M. V. Lawson and D. G. Jones; Graph inverse semigroups: Their characterization and completion; J. of Algebra, 409 (2014), 444-473

\Bibitem LV M. V.  Lawson and  A. Vdovina; The universal Boolean inverse semigroup presented by the abstract Cuntz-Krieger relations; arXiv:1902.02583v3
[math.OA], 2018

\Article MilanStein D. Milan, B. Steinberg; On inverse semigroup $C^*$-algebras and crossed products; Groups Geom. Dyn., 8 (2014), no. 2, 485-512

\Article StarlingTwo C. Starling; Boundary quotients of $\rm C^*$-algebras of right LCM semigroups; J. Funct. Anal., 268 (2015), no. 11, 3326-3356

\Article StarlingThree C. Starling; Inverse semigroups associated to subshifts; J. Algebra, 463 (2016), 211--233

\Article Stone M. H. Stone; Postulates for Boolean Algebras and Generalized Boolean Algebras; Amer. J. Math., 57 (1935), 703-732

\Bibitem Takesaki M. Takesaki; Theory of operator algebras. I; Springer-Verlag, New York-Heidelberg, 1979. vii+415 pp.

\Bibitem Wehrung F. Wehrung; Refinement monoids, equidecomposability types, and Boolean inverse semigroups; Lecture Notes in Mathematics, 2188, Springer, 2017

\end